\documentclass[regno, reqno,11pt]{amsart}
\usepackage{amsmath,amsthm,amssymb,amsfonts,epsfig,color,graphicx,enumerate,accents,subfigure,comment, esint}
\usepackage[a4paper,margin=2.59truecm]{geometry}
\usepackage{graphicx}
\usepackage{enumitem}

\usepackage[linkcolor=blue,colorlinks=true, citecolor = blue]{hyperref}

\DeclareMathOperator{\dive}{div}

\DeclareMathOperator{\dist}{dist}
\DeclareMathOperator{\supp}{supp}

\def\eps{{\varepsilon}}

\def\N{\mathbb{N}}
\def\NN{\mathcal{N}}
\def\R{\mathbb{R}}
\def\O{\Omega}

\def\vf{\varphi}
\def\AA{\mathcal{A}}
\def\HH{\mathcal{H}}

\def\DD{\mathcal{D}}

\def \< {\langle}
\def \>{\rangle}
\def\HH{\mathcal{H}}

\def\PP{\mathcal{P}}
\def\QQ{\mathcal{Q}}
\def\RR{\mathcal{R}}

\def\f{\phi}
\def\k{\kappa}

\newcommand{\be}{\begin{equation}}
\newcommand{\ee}{\end{equation}}

\newcommand{\bs}{\begin{split}}
\newcommand{\es}{\end{split}}

\newcommand{\Linfty}[2]{\| #1 \|_{L^{\infty}(#2)} }

\numberwithin{equation}{section}
\theoremstyle{plain}
\newtheorem{theorem}{Theorem}[section]
\newtheorem{lemma}[theorem]{Lemma}
\newtheorem{corollary}[theorem]{Corollary}
\newtheorem{proposition}[theorem]{Proposition}
\newtheorem{definition}[theorem]{Definition}
\theoremstyle{remark}
\newtheorem{remark}[theorem]{Remark}

\title[]{Tangential contact between free and fixed boundaries for variational solutions to variable coefficient Bernoulli type Free boundary problems}

\author{Diego Moreira$^{*}$}
\address{$^*$Departamento de Matemática, Universidade Federal da Ceara(Fortaleza, Brazil)}
\email{$^*$dmoreira@mat.ufc.br}
\author{Harish Shrivastava$^\dagger$}
\address{$^\dagger$Tata Institute of Fundamental Researcher-Centre of of Applicable Mathematics}
\email{$^\dagger$harish21@tifrbng.res.in, }
\begin{document}

\begin{abstract}
In this paper, we show that given appropriate boundary data, the free boundary and the fixed boundary of minimizers of functionals of type \eqref{functional} contact each other in a tangential fashion. We prove this result via classification of the global profiles, adapting the ideas from \cite{KKS06}.
\end{abstract}
\medskip

\maketitle


\textbf{2010 Mathematics Subjects Classification:}49J05, 35B65, 35Q92, 35Q35

\tableofcontents

\section{Introduction}
Objective of the paper of to study the behavior of free boundary near the fixed boundary of domain, for minimizers of Bernoulli type functionals with Hölder continuous coefficients. 
\be\label{functional}
J(v;A,\lambda_+,\lambda_-,\O)=\int_{\O} \Big ( \< A(x)\nabla v,\nabla v \>+ \Lambda(v)\Big )\,dx
\ee
$A$ is an elliptic matrix with Hölder continuous entries, and $\Lambda(v)=\lambda_+ \chi_{\{v>0\}}+\lambda_- \chi_{\{v\le0\}}$. We prove that if the value of boundary data and its derivative at a point are equal to zero (i.e. it satisfies the \eqref{DPT} condition mentioned below), then the contact of free boundary and the fixed boundary is tangential. 

Boundary interactions of free boundaries have gained significant attention in recent years. Whenever there are two medias involved, interactions of their respective diffusions  can be modeled by free boundary problems.  Often, free boundary of solution and fixed boundary of set come in contact. In applications, the Dam problem \cite{AG82} and Jets, Wakes and Cavities \cite{ZG57} model phenomenas which involve understanding of free boundary and fixed boundary.

Very recently, works of Indrei \cite{I00}, \cite{I19} study interactions of free boundaries and fixed boundaries for fully non-linear obstacle problems. We refer to \cite{GL19} where authors shed more light into angle of contact between fixed boundary and free boundary for one phase Bernoulli problem. 

As it is common by now, our strategy in this article is to classify blow up of minimizers by using the ideas from \cite{KKS06}. We prove that the blow up and also their positivity sets converge to a global solution in $\PP_{\infty}$ as defined in \cite{KKS06}. In Section \ref{setting up problem}, we list the assumptions and set some notations and then in Section \ref{blow up analysis}, we prove that blow ups of minimizers converge to that of global solutions (c.f. Definition 2.6). In the last section we prove our main result.

\section{Setting up the problem}\label{setting up problem}

We consider the following class of function which we denote as $\PP_r(\alpha, M,\lambda, \DD,\mu)$. Before definition, we set the following notation
\[
\begin{split}
B_R^+ := \Big \{ x\in B_R\mbox{ such that } x_N>0 \Big \}\\
B_R' := \Big \{ x\in B_R\mbox{ such that } x_N=0 \Big \}.
\end{split}
\]
For $x\in \R^N$ we denote $x'\in \R^{N-1}$ as the projection of $x$ on the plane $\{x_N=0\}$, we  denote the tangential gradient $\nabla'$ as follows
$$
\nabla' u := \Big ( \frac{\partial u}{\partial x_1}, ...\,\frac{\partial u}{\partial x_{N-1}} \Big ).
$$
 We define the affine space set $H^1_\f(B_R^+)$ as follows,
\be\label{Kf}
H^1_\f(B_R^+) = \left \{ v\in H^1(B_R^+) \,:\, v-\f\in H_0^1(B_R^+) \right \}.
\ee
For a given function $v\in H^1(B_2^+)$, we denote $F(v)$ as $F(v):=\partial \{v>0\}$ and $Id$ is the notation for $N\times N$ identity matrix. 

\begin{definition}
A function $u\in H^1(B_{2/r}^+)$ is said to belong to the class $\PP_r(\alpha, M,\lambda_{\pm}, \DD,\mu)$ if there exists $A\in C^{\alpha}(B_{2/r}^+)^{N\times N}$, $\f \in C^{1,\alpha }(B_{2/r}^+)$, $ \lambda_{\pm}>0$, $0<\mu<1$ and $\DD>0$ such that 
\begin{enumerate}[label=\textbf{$($P\arabic*$)$}]
\item \label{P1} $\Linfty{A}{B_{2/r}^+}\le M$, $\Linfty{\nabla \f}{B_{2/r}^+}\le M$, $[A]_{C^{\alpha}(B_{2/r}^+)},[\nabla \f]_{C^{1,\alpha}(B_{2/r}^+)}\le r^{\alpha}M$ \newline and $|\f(x')| \le Mr^{1+\alpha}|x'|^{1+\alpha}$ $($$x'\in B_{2/r}'$$)$. $\f$ satisfies the following Degenerate Phase Transition condition (\ref{DPT}) mentioned below.
\be\label{DPT}\tag{DPT}
\mbox{$\forall x'\in B_{2/r}'$ such that $\f (x')=0$, then $|\nabla' \f(x')|=0$}.
\ee
\item \label{P2} $A(0)= Id$, $\mu|\xi|^2 \le  \< A(x)\xi,\xi\> \le \frac{1}{\mu}|\xi|^2$ for all $x\in B_{2/r}^+$ and $\xi \in \R^N$.
\item \label{P3} $0<\lambda_-<\lambda_+$.
\item \label{P4} $u$ minimizes $J(\cdot; A,\lambda_+,\lambda_-,B_{2/r}^+)$ $($c.f. \eqref{functional}$)$ that is for every $u-v\in H_0^1(B_{2/r}^+)$
$$
\int_{B_{2/r}^+} \Big ( \< A(x) \nabla u,\nabla u  \>+ \Lambda (u) \Big ) \,dx \le \int_{ B_{2/r}^+} \Big ( \< A(x) \nabla v,\nabla v  \>+ \Lambda (v) \Big )\,dx \;\;
$$
$ (\Lambda(s)=\lambda_+ \chi_{\{s>0\}}+\lambda_- \chi_{\{s\le0\}}  )$ and $0\in F(u)\cap \overline{B_{2/r}^+}$.
\item \label{P5} $u \in H_{\f}^1(B_{2/r}^+)$. 
\item \label{P6} There exists $0<r_0$ such that for all $0<\rho\le r_0$ we have 
\be\label{density}
\frac{|B_\rho^+(0)\cap \{u>0\}|}{|B_\rho^+(0)|}>\DD.
\ee
\end{enumerate}
\end{definition}
\begin{remark}
In fact, the functions $u\in \PP_r(\alpha, M,\lambda_{\pm}, \DD,\mu)$ carry more regularity than being only a Sobolev function. They are Hölder continuous in $ (\overline{B_{2/r}^+})$ (c.f. Lemma \ref{lipschitz}). 
\end{remark}
In the absence of ambiguity on values of $\alpha, M,\lambda_{\pm}, \DD,\mu$ we use the notation $\PP_r$ in place of $\PP_r(\alpha, M,\lambda_{\pm}, \DD,\mu)$.  If $\f\in C^{1,\alpha}(B_2^+)$ satisfies (\ref{DPT}), from \cite[Lemma 10.1]{BM21}, we know that $\f^{\pm} \Big |_{B_2'} \in C^{1,\alpha}(B_2')$ and also
$$
\| \f^{\pm} \|_{C^{1,\alpha}(B_2')} \le \| \f \|_{C^{1,\alpha}(B_2')}.
$$
Given $v\in H^1(B_R^+)$ and $r>0$, we define the blow-up $v_r\in H^1(B_{R/r}^+)$ as follows  
\be\label{blowup}
v_r(x):=\frac{1}{r}v(rx).
\ee
For the coefficient matrix $A$, $A^r(x)$ is defined as follows
\be\label{Ar}
A^r(x):=A(rx).
\ee
One can check that if $u\in \PP_1(\alpha, M,\lambda_{\pm}, \DD,\mu)$, then $u_r \in \PP_r(\alpha, M,\lambda_{\pm}, \DD,\mu)$. Indeed if $u\in \PP_{1}$ and $u$ minimizer the functional $J$ (c.f. \ref{P4})
$$
J(v;A,\lambda_+,\lambda_-,B_2^+):= \int_{B_2^+}\Big ( \< A(x)\nabla v,\nabla v \> + \Lambda(v)\Big )\,dx,\qquad  (\Lambda(s)=\lambda_+ \chi_{\{s>0\}}+\lambda_- \chi_{\{s\le0\}}  )
$$
with boundary data $\f\in C^{1,\alpha}(B_2^+)$ (i.e. $u\in H_{\f}^1(B_2^+)$). Then by simple change of variables we can check that $u_r\in H^1_{\f_r}(B_{2/r}^+)$ (this verifies \ref{P5}) and $u_r$ minimizes 
$$
J(v;A^r,\lambda_+,\lambda_-,B_{2/r}^+):= \int_{B_{2/r}^+}\Big ( \< A^r(x)\nabla v,\nabla v \> + \Lambda(v)\Big )\,dx,\qquad  (\Lambda(s)=\lambda_+ \chi_{\{s>0\}}+\lambda_- \chi_{\{s\le0\}}  ).
$$
Moreover, if $A$ and $\f$ satisfy the conditions \ref{P1}, \ref{P2} for $r=1$, then $A^r$ and $\f_r$ satisfy \ref{P1}, \ref{P2} for $r$. \ref{P3} and \ref{P6} remains invariant under the change variables. Therefore $u_r\in \PP_r$.

In order to study the blow-up limits ($\lim_{r\to 0} u_r$) of functions $u\in \PP_1(\alpha, M,\lambda_{\pm}, \DD,\mu)$, we define a class of global solutions $\PP_{\infty}(C,\lambda_{\pm})$. Let us set the following notation before giving the definition
$$
\Pi:= \{x \in \R^N \,:\, x_N=0\}.
$$

\begin{definition}[Global solution]\label{global solution}
We say that $u\in H^1(\R^N_+)$ belongs to the class $\PP_{\infty}(C, \lambda_{\pm})$, that is, $u$ is a global solution if there exists $C>0$ and $0<\lambda_+<\lambda_-$ such that
\begin{enumerate}[label=\textbf{(G\arabic*)}]
\item \label{G1} $|u(x)|\le C|x|$, for all $x\in \R^N_+$,
\item \label{G2}$u$ is continuous up to the boundary $\Pi$,
\item \label{G3} $u=0$ on $\Pi$,
\item \label{G4} and for every ball $B_r(x_0)$, $u$ is a minimizer of $J(\cdot; {Id,\lambda_+,\lambda_-,B_r(x_0)\cap \R^N_+})$ $($c.f. \eqref{functional}$)$, that is 
$$
\int_{B_r(x_0)\cap \R^N_+} \Big ( |\nabla u|^2 + \Lambda(u)\Big )\,dx \le \int_{B_r(x_0)\cap \R^N_+} \Big ( |\nabla v|^2 +\Lambda(v) \Big )\,dx.
$$
Here $(\Lambda(s)=\lambda_+ \chi_{\{s>0\}}+\lambda_- \chi_{\{s\le0\}}  )$ and for every $v\in H^1(B_r(x_0)\cap \R^N_+)$ such that $u-v\in H_0^1(B_r(x_0)\cap \R^N_+)$.
\end{enumerate}
\end{definition}

Our main result intends to show that for a minimizer $u$ of $J(\cdot; A,\lambda_+,\lambda_-,B_2^+)$ with $A,\lambda_{\pm}$ and $u$ satisfying the properties \ref{P1}-\ref{P6}, the free boundary of every such minimizer touches the flat part of fixed boundary tangentially at the origin. For this, we prove that as we approach closer and closer to the origin, the free boundary points cannot lie outside any cone which is perpendicular to the flat boundary and has its tip at the origin. The main result in this paper is stated below.

\begin{theorem}\label{main result}
There exists a constant $\rho_0$ and a modulus of continuity $\sigma $  such that if 
$$u\in\PP_1(\alpha, M,\lambda_{\pm}, \DD,\mu)$$ 
then
$$
F(u)\cap B_{\rho_0}^+ \subset \{ x\,:\, x_N\le \sigma (|x|)|x| \}
$$
Here $\sigma$ depends only on $\alpha, M,\lambda_{\pm}, \DD,\mu$.
\end{theorem}

\section{Blow-up analysis}\label{blow up analysis}

The following is a classical result (c.f. \cite[Remark 4.2]{AC81}) , we present the proof for the case of variable coefficients.

\begin{lemma}\label{subharmonic}

Given a strictly elliptic matrix and bounded $A(x)$ and a non-negative continuous function $w$ such that $\dive(A(x)\nabla w)=0$ in $\{w>0\}\cap B_2^+$, then $w\in H^1_{loc}(B_2^+)$ and $\dive(A(x)\nabla w)\ge 0$ in weak sense in $B_2^+$.

\end{lemma}

\begin{proof}
Let $D\Subset B_2^+$ and $\eta\in C_c^{\infty}(B_2^+)$ be cutoff function for $D$. That is $\eta\in C_c^{\infty}(B_2^+)$ be such that 
$$
\eta(x) = \begin{cases}
1 \;\; \mbox{in $D$}\\
 0 \;\; \mbox {on $\partial B_2^+$}.
\end{cases}
$$
Since $\dive(A(x)\nabla w)=0$ in $\{w>0\}$, we have 
\[
\begin{split}
0&=\int_{B_2^+}\< A(x)\nabla w, \nabla \big (  (w-\eps)^+\eta^2  \big )\>\,dx\\
&=\int_{B_2^+\cap \{w>\eps\}} \eta^2 \< A(x) \nabla w, \nabla w \>\,dx+ \int_{B_2^+\cap \{w>\eps\}} w \< A(x) \nabla w, \nabla \eta^2 \>\,dx+\eps\int_{B_2^+ \cap \{w>\eps\}} \< A(x)\nabla w, \nabla \eta^2 \>\,dx
\end{split}
\]
which implies
\[
\begin{split}
  \int_{B_2^+ \cap \{w>\eps\}} \< A(x)\nabla w, \nabla w\> \eta^2 \,dx  \le  \int_{B_2^+\cap \{w>\eps\}} & \Big | w \< A(x) \nabla w, \nabla \eta^2 \>\Big | \,dx \\&+\eps\int_{B_2^+ \cap \{w>\eps\}} \Big | \< A(x)\nabla w, \nabla \eta^2 \>\Big | \,dx. 
\end{split}
\]
By the choice of $\eta$ and ellipticity of the matrix $A$, we obtain using Young's inequality 
\[
\begin{split}
\mu  \int_{B_2^+ \cap \{w>\eps\}}  |\nabla w|^2 \eta^2\,dx  &\le  \mu \Big | \int_{B_2^+ \cap \{w>\eps\}} \eta w |\nabla w| |\nabla \eta| \,dx \Big | + \eps \mu \Big | \int_{B_2^+ \cap \{w>\eps\}} \eta |\nabla w| |\nabla \eta|\,dx  \Big | \\
&\le C_1 (\mu) \Bigg [ \frac{1}{\delta}\int_{B_2^+ \cap \{w>\eps\}} w^2 |\nabla \eta|^2 \,dx + \delta \int_{B_2^+ \cap \{w>\eps\}} \eta^2 |\nabla w|^2 \,dx\\
&\qquad \qquad \qquad  \qquad +  \delta \eps \int_{B_2^+ \cap \{w>\eps\}} \eta^2 |\nabla w|^2\,dx +\frac{\eps}{\delta} \int_{B_2^+ \cap \{w>\eps\}}  |\nabla \eta|^2\,dx \Bigg]
\end{split}
\]
after choosing of $\delta>0$ very small and rearranging the terms in the equation above, since $\eta=1$ in $D$ we finally get 
\[
\begin{split}
\int_{D\cap \{w>\eps\}} |\nabla w|^2\,dx & \le \int_{B_2^+ \cap \{w>\eps\}}  |\nabla w|^2 \eta^2\,dx\\
& \le C(\mu) \Bigg [ \int_{B_2^+ \cap \{w>\eps\}} w^2 |\nabla n|^2 \,dx +  \int_{B_2^+ \cap \{w>\eps\}}  |\nabla \eta|^2\,dx \Bigg ]
\end{split}
\]
As $\eps\to 0$, we obtain 
\[
\begin{split}
\int_{D}|\nabla w|^2\,dx = \int_{\{w>0\} \cap D}|\nabla w|^2\,dx &= \lim_{\eps\to 0} \int_{\{w>\eps\}\cap D} |\nabla w|^2\,\\
&\le C(\mu) \lim_{\eps\to 0} \Bigg [  \int_{B_2^+ \cap \{w>\eps\}} w^2 |\nabla n|^2 \,dx +  \int_{B_2^+ \cap \{w>\eps\}}  |\nabla \eta|^2\,dx \Bigg ]\\
&\le C(\mu,D) \Bigg [ \int_{B_2^+\cap \supp(\eta)}w^2 \,dx + 1 \Bigg ].
\end{split} 
\]
Since $w\in C(B_2^+)$ therefore, $w$ is uniformly bounded in $\supp(\eta)$ and therefore 
$$
\int_{D}\Big ( |\nabla w|^2 + |w^2|  \Big )\,dx \le C(\mu,D) \Bigg [ \int_{B_2^+\cap \supp(\eta)}w^2 \,dx + 1 \Bigg ] \le C(\mu, D, \| w \|_{L^{\infty}(\supp(\eta))}).
$$
Now, for $0\le \varphi \in C_c^{\infty}(B_2^+)$ consider the test function 
$$
v = \varphi \Big ( 1 - \big (\min (2-\frac{w}{\eps} , 1) \big ) ^+    \Big ).
$$
$$
\int_{B_2^+}\< A(x) \nabla w, \nabla \varphi \>\,dx=\int_{B_2^+}\< A(x)\nabla w, \nabla \big  ( \varphi ((2-\frac{w}{\eps})\wedge 1  )^+ )\big )\>\,dx
$$
We can easily check that $v\ge 0$ in $B_2^+$ and $v\in H_0^1(B_2^+)$, in particular
$$
\varphi \Big ((2-\frac{w}{\eps})\wedge 1  \Big )^+ = 
\begin{cases}
\vf \;\;\;\qquad \qquad \;\; \qquad \qquad x\in \{w\le \eps\},\\
\vf \cdot \Big ( 2-\frac{w}{\eps} \Big )\; \qquad\;\;  \;\qquad x\in \{\eps<w\le 2\eps\},\\
0\;\;\;\; \qquad \qquad \qquad \qquad \;\;  x\in \{w>2\eps\}.
\end{cases}
$$
Therefore, we have
\[
\begin{split}
\int_{B_2^+}\< A(x) \nabla w, \nabla \varphi \>\,dx&=\int_{B_2^+}\< A(x)\nabla w, \nabla \big  ( \varphi ((2-\frac{w}{\eps})\wedge 1  )^+ )\big )\>\,dx\\
&= \int_{B_2^+\cap \{w\le \eps\}} \< A(x)\nabla w,\nabla \vf \>\,dx + \int_{B_2^+\cap \{\eps<w\le 2\eps\}} \< A(x)\nabla w,\nabla \Big ( \vf \big ( 2- \frac{w}{\eps}\big ) \Big ) \>\,dx\\
&= \int_{B_2^+\cap \{w\le \eps\}} \< A(x)\nabla w,\nabla \vf \>\,dx +2\int_{B_2^+\cap \{\eps<w\le 2\eps\}} \< A(x)\nabla w,\nabla \vf \>\,dx\\
&\qquad  \qquad \qquad \qquad \qquad \qquad \qquad \qquad-\frac{2}{\eps} \int_{B_2^+\cap \{\eps <w\le 2\eps\}} \< A(x)\nabla w,\nabla (w\vf) \>\,dx\\
&\le C(\mu) \int_{B_2^+ \cap \{\eps <u \le 2\eps\}} |\nabla w||\nabla \vf|\,dx+ \frac{2}{\eps} \int_{B_2^+ \cap \{\eps <w\le 2\eps\}}w\< A(x)\nabla w,\nabla \vf \>\,dx\\
&\qquad \qquad \qquad \qquad\qquad \qquad -\frac{2}{\eps} \int_{B_2^+ \cap \{ \eps <w\le 2\eps\}} \vf \< A(x)\nabla w,\nabla w\> \,dx\\
\le & C(\mu) \int_{B_2^+ \cap \{\eps<w\le 2\eps\}} |\nabla w||\nabla \vf|\,dx
\end{split}
\]
The last term goes to zero as $\eps\to 0$. Therefore, we can say that 
$$
\int_{B_2^+}\< A(x) \nabla w, \nabla \varphi \>\,dx\le 0.
$$
This concludes the proof.
\end{proof}

\begin{lemma}[Hölder continuity]\label{lipschitz}
If $u\in \PP_1$ then $u\in C^{0,\alpha_0}(\overline {B_2^+})$ for some $0<\alpha_0<1$. In fact,
$$
\|u\|_{C^{\alpha_0}(B_2^+)} \le C(\mu, \lambda_{\pm})\Linfty{u}{B_2^+}.
$$
\end{lemma}

\begin{proof}
The functional $J(\cdot;A, \lambda_+,\lambda_-,B_2^+)$ satisfy the hypothesis of \cite[Theorem 7.3]{eg05} and $\f\in C^{1,\alpha}(\overline{B_2^+})$, therefore Lemma \ref{lipschitz} follows from the arguments in \cite[Section 7.8]{eg05}.
\end{proof}

\begin{remark}
Since every function $u\in \PP_1$ is continuous. Therefore, the positivity set $\{u>0\}$ is an open set.
\end{remark}

\begin{corollary}\label{linear growth}
If $u\in \PP_1$, then $u^{\pm}$ are $A$-subharmonic.
\end{corollary}
\begin{proof}
The claim follows directly from Lemma \ref{lipschitz} and Lemma \ref{subharmonic}.
\end{proof}

\begin{lemma}\label{lemma 3.5}
If $u \in \PP_1$ then 
\be\label{lg}
|u(x)| \le C(\mu)M \,|x| \mbox{ in $B_1^+$.}
\ee.
\end{lemma}

\begin{proof}
Let $w$ be such that 
$$
\begin{cases}
\dive(A(x)w) = 0,\qquad \mbox{in $B_2^+$}\\
w = \f^+,\qquad \qquad\qquad\mbox{in $\partial B_2^+$}.
\end{cases}
$$
Since $u$ is $A$-subharmonic in $B_2^+$ (c.f. Corollary \ref{linear growth}), by maximum principle, if $x\in B_1^+$ we have
\be\label{borsuk 2}
\begin{split}
u^+(x)\le w(x)&\le w(x)-w(x')+w(x')\\
&\le \|\nabla w\|_{L^{\infty}(B_1^+)}|x-x'|+|\f^+(x')|\\
&=\|\nabla w\|_{L^{\infty}(B_1^+)}x_N+|\f^+(x')|\\
&\le (\|\nabla w\|_{L^{\infty}(B_1^+)} + M)|x|\;\forall x\in B_1^+
\end{split}
\ee
in the last inequality, we have used \ref{P1}. Now, we prove that the term $\|\nabla w\|_{L^{\infty}(B_1^+)}$ is uniformly bounded.

From \cite[Theorem 2]{B98}, we have 
\be \label{borsuk 1}
\|\nabla w\|_{L^{\infty}(B_1^+)}\le {C(\mu)} \Big [ \|w\|_{L^{\infty}(B_2^+)} + \|\f^+ \|_{C^{1,\alpha}(B_2')}  \Big ].
\ee
By comparison principle, $\| w\|_{L^{\infty}(B_2^+)} = \| w\|_{L^{\infty}(\partial B_2^+)}= \| \f^+ \|_{L^{\infty}(\partial B_2^+)}\le M$. Plugging this information in \eqref{borsuk 1}  
$$
\|\nabla w\|_{L^{\infty}(B_1^+)}\le C(\mu)M
$$
and then using \eqref{borsuk 2}, we obtain 
\be\label{u+}
u^+(x) \le C(\mu) M|x|\;\forall \;x\in B_1^+.
\ee
And analogously,  
\be\label{u-}
u^-(x) \le C(\mu) M|x|\;\forall \;x\in B_1^+.
\ee
By adding \eqref{u+} and \eqref{u-}, we prove \eqref{lg}.
\end{proof}

\begin{remark}\label{linear growth for blowup}
We can check that for every $u\in \PP_1$, $u_r\in \PP_r$ and $u_r$ is $A^r$-subharmonic and satisfies \eqref{lg} in $B_{1/r}^+$. That is 
$$
|u_r(x)| \le C(\mu)M\,|x|, \,\,x\in B_{1/r}^+.
$$
\end{remark}
\begin{lemma}[Uniform bounds in $H^1(B_1^+)$ norm]\label{global H1}\label{L2 bounds}
Let $u\in \PP_1$. Then for $R>0$ such that $2R\le \frac{2}{r}$, we have 
$$
\int_{B_{R}^+}|\nabla u_r|^2\,dx\le C (N, \lambda, \mu, R, M).
$$ 

\end{lemma}

\begin{proof}
Since $u_r\in \PP_r$, from \ref{P4} we can say that $u_r$ is a minimizer of $J(\cdot; A^r, \lambda_{\pm}, B_{2R}^+)$ with boundary data $\f_r$. Here $A^r$ and $\f_r$ satisfy the conditions \ref{P1} and \ref{P2}. Precisely speaking, $u_r$ is minimizer of the following functional
$$
J(v; A^r, \lambda_{\pm}, B_{2R}^+) : =\int_{B_{2R}^+} \Big ( \< A^r(x) \nabla v,\nabla v \> + \Lambda (v) \Big )\,dx  
$$
here $\big ( \Lambda(v) = \lambda_{+} \chi_{\{v>0\}} + \lambda_{-} \chi_{\{v\le0\}} \big )$. Consider $h\in H^1(B_{2R}^+))$ be a harmonic replacement 
$$
\begin{cases}
\dive(A^r(x)\nabla h) = 0\qquad \mbox{in $B_{2R}^+$}\\
h-u_r\in H_0^1(B_{2R}^+).
\end{cases}
$$
in other words, $h$ is the minimizer of $\int_{B_{2R}^+}\< A^r(x)\nabla h,\nabla h \> \,dx $, in the set $H_{u_r}^1(B_{2R}^+)$. 

From minimality of $u_r$ and the choice of $h$, we have
\[
\begin{split}
\int_{B_{2R}^+}\langle A^r(x) \nabla (u_r-h),\nabla &(u_r-h) \rangle\,dx = \int_{B_{2R}^+}\langle A^r(x) \nabla (u_r-h),\nabla (u_r+h-2h) \rangle\,dx\\
&= \int_{B_{2R}^+}\langle A^r(x) \nabla (u_r-h),\nabla (u_r+h) \rangle\,dx - 2\int_{B_{2R}^+}\langle A^r(x) \nabla (u_r-h),\nabla h \rangle\,dx\\
&= \int_{B_{2R}^+} \Big ( \langle A^r(x) \nabla u_r,\nabla u_r \rangle- \langle A^r(x) \nabla h,\nabla h \rangle \Big )\,dx\;\;\; \mbox{(since $h$ is $A^r$-harmonic in $B_R^+$)}\\
& \le \int_{B_{2R}^+} \Big ( \Lambda(h)-\Lambda(u_r) \Big )\,dx \le C(N,\lambda,R).
\end{split}
\]
We use ellipticity of $A$, we get 
$$
 \int_{B_{ R}^+}|\nabla (u_r-h)|^2\,dx \le \int_{B_{2R}^+}|\nabla (u_r-h)|^2\,dx\le C(N, \lambda, \mu, R)
$$
expanding the left hand side, we get 
\[
\begin{split}
 \int_{B_{R}^+}|\nabla u_r|^2\,dx\le  \int_{B_{R}^+}|\nabla u_r|^2\,dx+\int_{B_{R}^+}|\nabla h|^2\,dx&\le C(N, \lambda, \mu, R)+2\int_{B_{R}^+}\nabla u_r\cdot \nabla h\,dx\\
 &\le C(N, \lambda, \mu, R) + \eps  \int_{B_{R}^+}|\nabla u_r|^2\,dx +\frac{1}{\eps} \int_{B_{R}^+}|\nabla h|^2\,dx
\end{split}
\]
by choosing $\eps=\frac{1}{8}$ we are left with the following,
$$
 \int_{B_{R}^+}|\nabla u_r|^2\,dx\le C(N, \lambda, \mu, R) \Big ( 1 +   \int_{B_{R}^+}|\nabla h|^2\,dx \Big ).
$$
From \cite[Theorem 2]{B98}, $\Linfty{\nabla h}{B_{R}^+}\le C(\mu,M)$. Thus we obtain a uniform bound on $\int_{B_{R}^+}|\nabla h|^2\,dx$.
\end{proof}

\begin{lemma}[Compactness]\label{convergence}
Let $r_j\to 0^+$, and a sequence $\{v_j\} \in \PP_1$. Then the blow-ups $u_j : = (v_j)_{r_j}$ $($as defined in \eqref{blowup}$)$ converges $($up to subsequece$)$ uniformly in $B_R^+$  and weakly in $H^1(B_R^+)$ to some limit for any $R>0$. Moreover, if $u_0$ is such a limit of $u_j$ in the above mentioned topologies, then $u_0$ belongs to $\PP_{\infty}$.
\end{lemma}
\begin{proof}
We fix $R>0$, since $v_j\in \PP_1$, therefore $u_j\in \PP_{r_j}$ and as argued in the proof of previous Lemma, the functions $u_j$ are minimizers of the functional $J (\cdot; A_j, \lambda_{\pm}, B_{R}^+)$ for $j$ sufficiently large that $R<\frac{1}{r_j}$. We set the notation for the functional $J_j$ as
\be\label{Jj}
J_j(v):=J (v; A_j, \lambda_{\pm}, B_{R}^+) = \int_{B_{R}^+} \Big ( \langle A_{j}(x) \nabla v,\nabla v\rangle+\Lambda(v) \Big )\,dx,\qquad v\in H_{u_j}^1(B_{R}^+).
\ee
We also denote the boundary values for $u_j \in \PP_{r_j}$ as $\f_j$. Here the sequences $A_j\in {C^{\alpha}(B_{2/r_j}^+)}^{N\times N}$ and $\f_j\in C^{1,\alpha}(B_{2/r_j})$ satisfy the condition \ref{P1}, \ref{P2} with $r=r_j$, $\Lambda(v) := \lambda_{+} \chi_{\{v>0\}} + \lambda_{-} \chi_{\{v\le0\}} $. We set the following notation for the functional $J_0$
\be\label{J0}
J_0(v;B_R^+):=\int_{B_{R}^+}|\nabla v|^2+\Lambda(v)\,dx.
\ee

From Lemma \ref{lipschitz}, we know that $u_j\in C^{\alpha_0}(\overline{B_{2/r_j}^+)}$ which implies $C^{\alpha_0}(\overline{B_{R}^+)}$. In particular $\|u_j\|_{C^{0,\alpha_0}(B_R^+)} \le C(\mu,\lambda_{\pm})$. Hence, $u_j$ is a uniformly bounded and equicontinuous sequence in $\overline{B_R^+}$, we can apply Arzela Ascoli theorem to show that $u_j $  uniformly converges to a function $u_0 \in C^{0,\alpha_0}(\overline{B_R^+})$.

Since $u_j=\f_j$ on $B_{R}'$, from \ref{P1} we  have $|\f_j(x)| \le M r_j^{1+\alpha} |x|^{1+\alpha}$ for $x\in B_{R}'$, therefore \newline $|\f_j(x)| \le C(M,\alpha) r_j^{1+\alpha} R^{1+\alpha}$. Hence $\f_j\to 0$ uniformly on $B_R'$. We have 
$$ u_0 =\lim_{j\to \infty} u_j = \lim_{j\to \infty} \f_j =0 \mbox{ on $B_R'$.}$$
 Thus $u_0$ satisfies \ref{G2} and \ref{G3} inside the domain $\overline{B_R^+}$. Also, from Lemma \ref{L2 bounds} we have 
\be\label{bound}
\int_{B_R^+} |\nabla u_j|^2\,dx  \le C (N, \lambda_{\pm}, \mu, R, M, \alpha)
\ee
then, by the linear growth condition (c.f. Remark \ref{linear growth for blowup}), $u_j$ also satisfies 
\be
|u_j(x)| \le C(\mu,\alpha)M|x|\;\;x\in B_R^+.
\ee
Hence, passing to the limit, we have $|u_0(x)| \le C(\mu,\alpha)M|x|,\;\forall x\in B_R^+$. In other words $u_0$, satisfies \ref{G1} in $\overline{B_R^+}$. Moreover, we have 
\be\label{bound L2}
\int_{B_R^+}|u_j|^2\,dx\le C(\mu,\alpha,M) \int_{B_R^+} |x|^2\,dx  \le C(\mu, \alpha, M, N, R).
\ee
Thus \eqref{bound} and \eqref{bound L2} imply that $u_j$ is a bounded sequence in $H^1(B_R^+)$. Hence, up to a subsequence, $u_j \rightharpoonup u_0$ weakly in $H^1(B_R^+)$. We rename the subsequence again as $u_j$.

We have found a blow-up limit up to a subsequence $u_0$ and have shown that $u_0$ satisfies \ref{G1}, \ref{G2} and \ref{G3} in $\overline{B_R^+}$. In order to show that $u_0\in \PP_{\infty}$, it only remains to verify that $u_0$ satisfies \ref{G4}, i.e. $u_0$ is a local minimizer of $J_{0}(\cdot; B_R^+)$ for all $R>0$ (c.f. \eqref{J0}). For that, we first claim that 
\be\label{gamma convergence}
\int_{B_R^+} \Big ( |\nabla u_0|^2+\Lambda(u_0) \Big )\,dx  \le \liminf_{j\to \infty} \int_{B_R^+} \Big ( \langle A_{j}(x) \nabla u_j,\nabla u_j\rangle+\Lambda(u_j) \Big )\,dx.
\ee
Indeed, let us look separately at the term $J_j(u_j)$ on the right hand side of the above equation

$$
J_{j}(u_j)= \int_{B_R^+} \Big ( \< A_j(x) \nabla u_j,\nabla u_j \>  + \lambda_+ \chi_{\{u_j>0\}}  + \lambda_- \chi_{\{u_j\le 0\}}\Big ) \,dx.
$$

We rewrite the first term as follows 
\be\label{difference}
\int_{B_R^+}  \< A_j(x) \nabla u_j,\nabla u_j \> \,dx = \int_{B_R^+} \< (A_j(x) -Id )\nabla u_j,\nabla u_j \> \,dx+\int_{B_R^+}|\nabla u_j|^2\,dx.
\ee
From \ref{P1} and \ref{P2} we have for all $x\in B_R^+$
$$
\|A_j(x) - Id\|_{L^{\infty}(B_R^+)} \le M r_j^{\alpha} |x|^{\alpha}\le C(M,R,\alpha)r_j^{\alpha} \to 0\mbox{ as $j\to \infty$}.  
$$
Therefore $A_j\to Id$ uniformly and $\|\nabla u_j\|_{L^2(B_2^+)}$ is bounded (c.f. \eqref{bound}). Hence, the first term on the right hand side of \eqref{difference} tends to zero as $j\to \infty$. Thus, from \eqref{difference} and by weak lower semi-continuity of $H^1$ norm, we have
\be\label{eq10}
\int_{B_R^+}|\nabla u_0|^2\,dx
\le \liminf_{j\to \infty} \int_{B_R^+} |\nabla u_j|^2\,dx =  \liminf_{j\to \infty} \int_{B_R^+} \< A_j(x) \nabla u_j,\nabla u_j \> \,dx  .
\ee
For the second term, we claim that
\be\label{last term}
\begin{split}
\int_{B_R^+}\lambda_+\chi_{\{u_0>0\}}+ \lambda_-\chi_{\{u_0\le0\}} \,dx &\le \liminf_{j\to \infty} \int_{B_R^+} \lambda_+\chi_{\{u_j>0\}}+\lambda_-\chi_{\{u_j\le0\}}\,dx.\\
\end{split}
\ee 
To see this, we first show that for almost every $x\in B_R^+$, we have
\be\label{char}
\lambda_{+}\chi_{\{u_0>0\}}(x)+ \lambda_- \chi_{\{u_0\le0\}}(x) \le \liminf_{j\to \infty} \left ( \lambda_+ \chi_{\{u_j>0\}}(x)+\lambda_-\chi_{\{u_j\le0\}}(x)\right ).
\ee
Indeed, let $x_0\in B_R^+\cap \big ( \{u_0>0\}\cup \{u_0<0\} \big )$. Then by the uniform convergence of $u_j$ to $ u_0$, we can easily see that $u_j(x_0)$ attains the sign of $u_0(x_0)$ for sufficiently large value of $j$. Hence, \eqref{char} holds in $\{u_0>0 \} \cup \{u_0<0\}$.

Now, assume $x_0\in B_R^+\cap \{u_0=0\}$. Then left hand side of \eqref{char} is equal to 
$$
\lambda_{+}\chi_{\{u_0>0\}}(x_0)+ \lambda_- \chi_{\{u_0\le0\}}(x_0) = \lambda_{-}.
$$
Regarding RHS of \eqref{char}, we see that 
$$
\lambda_+ \chi_{\{u_j>0\}}(x_0)+\lambda_-\chi_{\{u_j\le0\}}(x_0) = 
\begin{cases}
 \lambda_+, \qquad \mbox{if $u_j(x_0)>0$}\\
  \lambda_-, \qquad \mbox{if $u_j(x_0)\le 0$}.
\end{cases}
$$
Since $\lambda_- <\lambda_+$ (c.f. \ref{P3}), the right hand side in the equation above is always greater than or equal to $ \lambda_- $. Then 
\[
\begin{split}
\lambda_{+}\chi_{\{u_0>0\}}(x_0)+ \lambda_- \chi_{\{u_0\le0\}}(x_0) = \lambda_{-} \le \liminf_{j\to \infty} \left ( \lambda_+ \chi_{\{u_j>0\}}(x_0)+\lambda_-\chi_{\{u_j\le0\}}(x_0) \right ).
\end{split}
\]
Thus, \eqref{char} is proven for all $x\in B_R^+$ and hence \eqref{last term} holds by Fotou's lemma. 

By adding  \eqref{eq10} and \eqref{last term} and \cite[Theorem 3.127]{D17} we obtain \eqref{gamma convergence}.  Now we will use \eqref{gamma convergence} prove the minimality of $u_0$ for the functional $J_0(\cdot;B_R^+)$ (c.f. \ref{J0}).

Pick any $w\in H^1(B_R^+)$ such that, $u_0-w\in H_0^1(B_R^+)$.  We construct an admissible competitor $w_j^{\delta}$ to compare the minimality of $u_j$ for the functional $J_j(\cdot;B_R^+)$. Then we intend to use \eqref{gamma convergence}.

In this direction, we define two cutoff functions $\eta_{\delta}: \R^N\to \R$ and $\theta :\R \to \R$ as follows,
$$
\eta_{\delta}(x):=
\begin{cases}
1,\; x\in B_{R-\delta}\\
0, \;x \in \R^N\setminus B_{R}
\end{cases},
\theta(t):=
\begin{cases}
1,\;|t|\le 1/2\\
0,\; |t| \ge 1.
\end{cases}
$$
we can take $|\nabla \eta_{\delta}|\le \frac{C(N)}{\delta}$.  We define $\theta_j(x)= \theta(\frac{x_N}{d_j})$, for a sequence $d_j\to 0$, which we be suitably chosen in later steps of the proof. Let $w_j^{\delta}$ be a test function defined as 
\be\label{test function}
w_j^{\delta} := w + (1-\eta_{\delta}) (u_j-u_0) +\eta_{\delta} \theta_{j} \f_j.
\ee
Since, the function $w_j^{\delta} -w = (1-\eta_{\delta}) (u_j-u_0)+\eta_{\delta} \theta_{j} \f_j$ is continuous in $\overline{B_R^+}$ and is pointwise equal to zero on $\partial B_R^+$, which is a Lipschitz surface in $\R^N$. Therefore, $u_j-w_j^{\delta} \in H_0^1(B_R^+)$. For further steps, the reader can refer to the Figure 1.
 
 \begin{figure}[h]\label{fig1}
\centering
\includegraphics[width=0.7\textwidth]{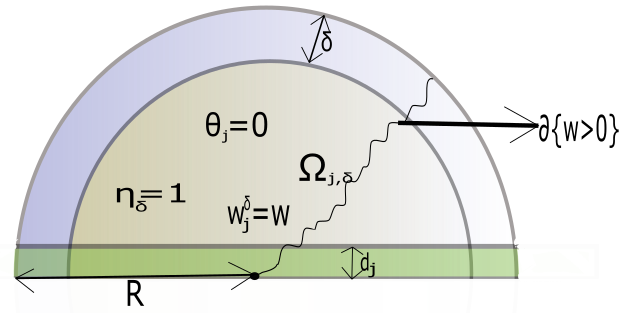}
\caption{(curvy line represents the free boundary of $w$)}
\end{figure}

Let $\O_{\delta,j}= B_R^+ \cap \{ \theta_j=0 \} \cap \{\eta_{\delta}=1\}$,  and $\RR_{\delta,j} = B_R^+ \setminus \O_{\delta,j}$ by observing $w_j^{\delta} = w$ on $\O_{\delta,j}$ we see that 
\[
\begin{split}
|\{w_j^{\delta}>0 \} \cap B_R^+| &= |\{w_j^{\delta}>0 \}\cap \O_{\delta,j} |+ |\{w_j^{\delta}>0 \}\cap \RR_{\delta,j}|\\
&= |\{w>0 \}\cap \O_{\delta,j} |+ |\{w_j^{\delta}>0 \}\cap \RR_{\delta,j}|\\
&= |\{w>0 \}\cap (B_R^+\setminus \RR_{\delta,j} )|+ |\{w_j^{\delta}>0 \}\cap \RR_{\delta,j}|\\
&= |\{w>0 \}\cap B_R^+| - |\{w>0 \} \cap \RR_{\delta,j}| + |\{w_j^{\delta}>0 \}\cap \RR_{\delta,j}|.
\end{split}
\]
From the above discussions, we have
\[
\begin{split}
|\{w>0\}\cap B_R^+| - |\RR_{\delta,j}| \le |\{w_j^{\delta}>0\}\cap B_R^+| \le |\{w>0\}\cap B_R^+| + |\RR_{\delta,j}|.
\end{split}
\]
Since we know that $\lim_{\delta\to 0}\left ( \lim_{j\to 0} |\RR_{\delta,j}| \right )=0$, we have 
\be\label{eq15}
\lim_{\delta \to 0}\left ( \lim_{j\to \infty} |\{ w_{j}^{\delta}>0 \} \cap B_R^+| \right )= |\{w>0\} \cap B_R^+|
\ee
and similarly
\be\label{eq16}
\lim_{\delta \to 0}\left ( \lim_{j\to \infty} |\{ w_{j}^{\delta}\le 0 \} \cap B_R^+| \right )= |\{w\le 0\} \cap B_R^+|.
\ee
Given $u_j\in \PP_{r_j}$ and $w_j^{\delta}-u_j \in H_0^1(B_R^+)$, from the minimility of $u_j$ for the functional $J_j$ we  have 
$$
\int_{B_R^+} \Big ( \<  A_j(x) \nabla u_j, \nabla u_j \> +\lambda_+ \chi_{\{u_j>0\}} + \lambda_- \chi_{\{u_j\le 0\}} \Big ) \,dx \le  \int_{B_R^+} \Big ( \<  A_j(x) \nabla w_j^{\delta}, \nabla w_j^{\delta} \> +\lambda_+ \chi_{\{w_j^{\delta}>0\}} + \lambda_- \chi_{\{w_j^{\delta}\le 0\}} \Big ) \,dx
$$
and from \eqref{gamma convergence} we obtain 
\be\label{eq12}
\begin{split}
\int_{B_R^+} \Big ( |\nabla u_0|^2 + \lambda_+ \chi_{\{u_0>0\}}+\lambda_- \chi_{\{u_0\le 0\}}\Big ) \,dx & \le \liminf_{j\to \infty}  \int_{B_R^+}\Big ( \<  A_j(x) \nabla w_j^{\delta}, \nabla w_j^{\delta} \> +\lambda_+ \chi_{\{w_j^{\delta}>0\}} + \lambda_- \chi_{\{w_j^{\delta}\le 0\}}\Big ) \,dx\\
&\le \limsup_{j\to \infty}  \int_{B_R^+}\Big ( \<  A_j(x) \nabla w_j^{\delta}, \nabla w_j^{\delta} \> +\lambda_+ \chi_{\{w_j^{\delta}>0\}} + \lambda_- \chi_{\{w_j^{\delta}\le 0\}}\Big ) \,dx
\end{split}
\ee
from the same reasoning as for the justification of \eqref{eq10}, we have 
\be\label{eq13.0}
 \limsup _{j\to \infty} \int_{B_R^+} \<  A_j(x) \nabla w_j^{\delta}, \nabla w_j^{\delta} \> \,dx = \limsup_{j\to \infty} \int_{B_R^+}|\nabla w_j^{\delta}|^2\,dx
\ee
therefore rewriting \eqref{eq12}
\be\label{eq13}
\int_{B_R^+} \Big ( |\nabla u_0|^2 + \lambda_+ \chi_{\{u_0>0\}}+\lambda_- \chi_{\{u_0\le 0\}} \Big ) \,dx \le  \limsup_{j\to \infty} \int_{B_R^+} \Big (  |\nabla w_j^{\delta}|^2\,dx+\lambda_+ \chi_{\{w_j^{\delta}>0\}} + \lambda_- \chi_{\{w_j^{\delta}\le 0\}} \Big ) \,dx.
\ee
We claim that
\be\label{eq14}
\lim_{\delta \to 0}\left ( \limsup_{j\to \infty} \int_{B_R^+} | \nabla w_j^{\delta}|^2 \,dx  \right )=  \int_{B_R^+} |\nabla w|^2\,dx
\ee
To obtain the claim above, we prove that 
$$
\lim_{\delta \to 0}\left (  \limsup_{j\to \infty} \int_{B_R^+} |\nabla (w_j^{\delta}-w)|^2\,dx \right ) = 0.
$$
From the definition of $w_j^{\delta}$, we know that  
$$
w^{\delta}_j-w = (1-\eta_{\delta}) (u_j-u_0) +\eta_{\delta} \theta_{j} \f_j
$$
therefore we have 
\be\label{all terms}
\begin{split}
\int_{B_R^+} |\nabla (w_j^{\delta}-w)|^2\,dx & \le C \Bigg (  \int_{B_R^+} |\nabla  \big ( (1-\eta_{\delta})(u_j-u_0) \big )|^2 \,dx + \int_{B_R^+}  |\nabla (\theta_j \eta_{\delta} \f_j)|^2\,dx \Bigg )\\
&\le  C(N) \Bigg ( \int_{B_R^+} (1-\eta_{\delta})^2 |\nabla (u_j-u_0)|^2\,dx+\frac{1}{\delta^2}\int_{B_R^+}|u_j-u_0|^2\,dx\\
& \qquad  \qquad \qquad \qquad \qquad \qquad \qquad \qquad \qquad+ \int_{B_R^+}  |\nabla (\theta_j \eta_{\delta} \f_j)|^2\,dx \Bigg ).
\end{split}
\ee
Let us consider the first term on the right hand side. We know that $\int_{B_R^+} |\nabla (u_j-u_0)|^2\,dx$ is uniformly bounded in $j\in \N$ (c.f. \eqref{bound}). Therefore 
\be\label{t1}
\begin{split}
 \lim_{\delta\to 0 } \Big ( \limsup_{j\to \infty} \int_{B_R^+} (1-\eta_{\delta})^2 |\nabla (u_j-u_0)|^2\,dx \Big ) &= \Big (\lim_{\delta\to 0}\Linfty{1-\eta_{\delta}}{B_R^+}^2 \Big )  \Big ( \limsup_{j\to \infty} \int_{B_R^+} |\nabla (u_j-u_0)|^2\,dx \Big )\\
&\le C(N,\lambda,\mu,R,M,\alpha)\lim_{\delta\to 0}\Linfty{1-\eta_{\delta}}{B_R^+} =0
\end{split}
\ee
Regarding the second term, since $|u_j-u_0|$ tends to zero in $L^2(B_R^+)$ as $j\to \infty$, therefore the second term also tends to zero as $j\to \infty$. We write 
\be\label{t2}
\lim_{\delta\to 0 } \Big ( \lim_{j\to \infty}\frac{1}{\delta^2}\int_{B_R^+}|u_j-u_0|^2\,dx \Big ) =0.
\ee
Lastly, we claim that 
\be\label{t3}
\lim_{j\to \infty}\int_{B_R^+}  |\nabla (\theta_j \eta_{\delta} \f_j)|^2\,dx =0.
\ee
Indeed, we have
\be\label{final part}
\int_{B_R^+}  |\nabla (\theta_j \eta_{\delta} \f_j)|^2\,dx\le C \left (\int_{B_R^+} |\nabla \eta_{\delta}|^2 (\theta_j \f_j)^2\,dx +\int_{B_R^+}|\nabla \f_j|^2 (\eta_{\delta} \theta_j)^2 \,dx+\int_{B_R^+}  |\nabla \theta_j |^2(\eta_{\delta} \f_j)^2\,dx \right ).
\ee
Since $\eta_{\delta}, \theta_j\le 1$, $|\nabla \eta_{\delta}| \le \frac{C(N)}{\delta}$ and $\Linfty{\nabla \f_j}{B_R^+}\le M$ (c.f. \ref{P1}), we obtain  
\be\label{dterm}
\begin{split}
\int_{B_R^+} |\nabla \eta_{\delta}|^2 (\theta_j \f_j)^2\,dx +\int_{B_R^+}|\nabla \f_j|^2 (\eta_{\delta} \theta_j)^2 \,dx &\le C(N) |\{ \theta_j \neq 0 \}\cap B_R^+| \Big ( \frac{1}{\delta^2} + M^2 \Big )\\
\end{split}
\ee
We know that $|\{ \theta_j \neq 0\}\cap B_R^+|\to 0$ as $j\to \infty$, hence from \eqref{dterm}, the first and second term in \eqref{final part} tend to zero as $j\to \infty$. The last term in \eqref{all terms} also tends to zero as $j\to \infty$, indeed from \ref{P1} we have $[\nabla \f_j]_{C^{\alpha}(B_R^+)} \le r_j^\alpha M$. Since $\f_j(0) = 0$, therefore we have $|\f_j|\le R^{\alpha} [\nabla \f_j]_{B_R^+}\le R^{\alpha} r_j^{\alpha} M$ in $B_R^+$. Also, observing that $|\nabla \theta_j| \le \frac{1}{d_j}$, $\eta_{\delta}\le 1$ in $B_{R}^+$ we have
$$
\int_{B_R^+}  |\nabla \theta_j |^2 (\eta_{\delta} \f_j)^2\,dx\le M R^{2\alpha}\frac{r_j^{2 \alpha}}{d_j^2} |B_R^+| 
$$
if we choose a sequence $d_j\to 0^+$ such that we also have $\frac{r_j^{\alpha}}{d_j} \to 0$, the third term in \eqref{final part} tends to zero as $j\to \infty$. Plugging in the estimates above \eqref{t1}, \eqref{t2}, \eqref{t3} in \eqref{all terms}, we obtain the claim (\ref{eq14}). 

From the equations (\ref{eq15}), (\ref{eq16}), and (\ref{eq14}) we obtain that the right hand side of \eqref{eq13} is equal to $J_0(w;B_R^+)$, therefore $u_0$ is a minimizer of $J_0(\cdot;B_R^+)$. That is 
$$
J_0(u_0;B_R^+) \le J_0(w;B_R^+)\;
$$
for every $w \in H^1(B_R^+)$ such that $u_0-w\in H_0^1(B_R^+)$. Since the inequality above (which corresponds to \ref{G4}) and other verified properties of $u_0$ (i.e. \ref{G1}, \ref{G2} and \ref{G3} in $B_R^+$) hold for every $R>0$, therefore $u_0 \in H_{loc}^1(\R^N_+)$ satisfies all the properties in the Definition \ref{global solution}. Hence $u_0\in \PP_{\infty}$.
 
\end{proof}

After proving that the (subsequential) limits of blow-up are global solutions, we proceed to show that the positivity sets (and hence the free boundaries) of blow-ups converge in certain sense to that of blow-up limit. For this we will need to establish that the minimizers $u\in \PP_r$ are non-degenerate near the free boundary. In the proof below, we adapt the ideas from \cite{ACF84}.

\begin{proposition}[Non-degeneracy near the free boundary] \label{non degeneracy}
Let $u\in \PP_{r_0}$ for some $r_0>0$ and $x_0\in B_{2/r_0}^+$. Then, for every $0<\kappa<1$ there exists a constant $c(\mu,N,\k,\lambda_{\pm})>0$ such that for all $B_{r}(x_0)\subset B_{\frac{2}{r_0}}^+$, we have 
\be
\frac{1}{r} \fint _{\partial B_r(x_0)} u^+\,d\HH^{N-1}(x)\,<\,c(\mu,N,\k,\lambda_{\pm}) \implies \mbox{$u^+=0$ in $B_{\kappa r}(x_0)$}.
\ee 
\end{proposition}

\begin{proof}
We fix $x_0\in \{u>0\}\cap B_2^+$ and $r>0$ such that $B_r(x_0) \subset B_2^+$. We denote \newline $\gamma := \frac{1}{r} \fint _{B_r(x_0)} u^+\,dx$. We know from Lemma \ref{lipschitz} that the set $\{u>0\}$ is open. Also, since $u\in \PP_{r_0}$, there exists $A\in C^{\alpha}(B_{2/r_0}^+)^{N\times N}$, $\vf \in C^{1,\alpha}(B_{2/r_0}^+)$, $\lambda_{\pm}$ satisfying \ref{P1}-\ref{P6}. Therefore $u$ solves the PDE $\dive(A(x)\nabla u) = 0$ in $\{u>0\}\cap B_{2/r_0}^+$. By elliptic regularity theory, $u$ is locally $C_{loc}^{1,\alpha}(\{u>0\}\cap B_{2/r_0}^+)$. Then, for almost every $\epsilon>0$, $B_r\cap \partial \{u>\eps\}$ is a $C^{1,\alpha}$ surface. Pick one such small $\eps>0$ and we consider the test function $v_{\eps}$ given by
$$
\begin{cases}
\dive(A(x)\nabla v_{\eps}) =0 \qquad \mbox{in $(B_r(x_0)\setminus B_{\kappa r}(x_0))\cap \{u>\eps\}$}\\
v_{\eps} = u\qquad \qquad \;\; \qquad \qquad \mbox{in $B_{r}(x_0)\cap \{u\le \eps \}$}\\
v_{\eps}=\eps \qquad \;\; \qquad \qquad \qquad\mbox{in $B_{\k r}(x_0)\cap \{u>\eps\}$}\\
v_{\eps}=u \qquad \; \; \qquad \qquad \qquad \mbox{on $\partial B_r(x_0)$.}
\end{cases}
$$ 
The function $v_{\eps}$ defined above belongs to $H^1(B_r(x_0))$, thanks to \cite[Theorem 3.44]{DD12} $\big ($\cite[Theorem 3.44]{DD12} is proven for $C^1$ domains, but the proof can also be adapted for Lipschitz domains  \cite[Theorem 4.6]{GE15}$\big )$.
 \begin{figure}[h]\label{fig1}
\centering
\includegraphics[width=0.7\textwidth]{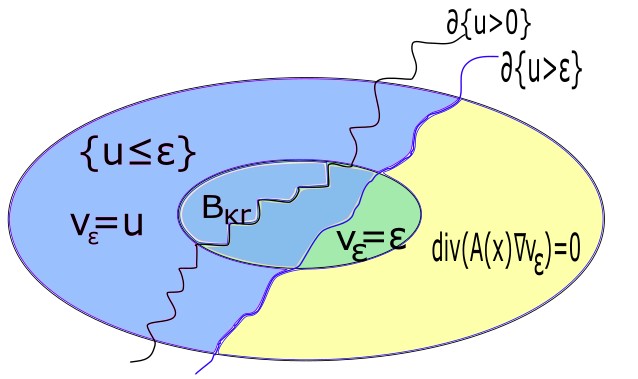}
\caption{Graph of $v_{\eps}$.}
\end{figure}
We intend to show that $ v_{\eps}$ is bounded in $H^{1}(B_r(x_0))$. This ensures the existence of limit $\lim_{\eps\to 0}v_{\eps}$ exists in weak sense in $H^1(B_r(x_0))$ and strong sense in $L^2(B_r(x_0))$. Let $G$ be the Green function for $L(v)=\dive(A(x)\nabla v)$ in the ring $B_r(x_0)\setminus B_{\k r}(x_0)$. Then if there is a function $w$ such that 
\be\label{PDE w}
\begin{cases}
\dive(A(x)\nabla w)=0 \; &\mbox{in $B_r(x_0)\setminus B_{\k r}(x_0)$}\\
w=u \; &\mbox{on $\partial B_r(x_0)\cap \{u>\eps\}$}\\
w=\eps \; &\mbox{elsewhere on $\partial (B_r(x_0)\setminus B_{\k r}(x_0)).$}
\end{cases}
\ee
We can also write that $w-\eps = (u-\eps)^+$ on $\partial B_r(x_0)$ and $w-\eps =0$ on $\partial B_{\k r}(x_0)$. Consider any sequence $\{x_k\} \subset B_r(x_0)\setminus B_{\k r}(x_0)$ such that $x_k \to \bar x \in \partial B_{\k r}(x_0)$.
By Green representation formulae for $w-\eps$ in \eqref{PDE w}, we have $w(x_k) \to w(\bar x)$, indeed since  
\[
\begin{split}
0=(w-\eps)(\bar x)&=\lim_{k\to \infty}(w-\eps)(x_k)\\
&=\lim_{k\to \infty}\int_{\partial (B_r(x_0)\setminus B_{\k r}(x_0))}\left(u-\eps)^+(y) (A(y)\nabla_yG(x_k,y)\right )\cdot \nu_y\,d \sigma(y)\\
&=   \lim_{k\to \infty}\int_{\partial B_r(x_0) \cap \{u>\eps\}}\left(u-\eps)^+(y) (A(y)\nabla_yG(x_k,y)\right )\cdot \nu_y\,d \sigma(y)  ,\qquad \forall x\in \partial B_{\k r}(x_0)
\end{split}
\]
where $\nu_y$ is the unit outer normal vector at a point $y$ on the boundary. 
We apply same arguments as above to $\nabla w(\bar x)$ and from \cite[Theorem 3.3 (vi)]{GW82} on $G(x,y)$ and therefore for $\bar x\in \partial B_{\k r}(x_0)$
\be\label{green estimates}
\begin{split}
\big | \nabla w(\bar x)\big | &\le C(\mu) \lim_{k\to \infty} \int_{\partial B_r(x_0)\cap \{u>\eps\}}\Big | \nabla_x \left ( \frac{\partial}{\partial _{\nu_y}}G(x_k,y) \right ) (u-\eps)^+ \Big |\,dx\\
& \le C(\mu) \lim_{k\to \infty} \int_{\partial B_r(x_0)} \frac{1}{|x_k-y|^N}(u-\eps)^+\,d\HH^{N-1}(y)\\
&\le \frac{C(\mu,N)}{(1-\k)^N}\frac{1}{r}\fint _{\partial B_r} (u-\eps)^+\,d\HH^{N-1}(y)\le C(\mu,N,\k) \gamma\;\mbox{on $\partial B_{\k r}(x_0)$}.
\end{split}
\ee
We can easily check by respective definitions that $w\ge v_{\eps}$ on $\partial (B_r(x_0)\setminus B_{\k r}(x_0))$, moreover, by maximum principle, since $\dive(A(x)\nabla w) =0$ in $B_r(x_0)\setminus B_{\k r}(x_0)$ and $w\ge \eps$ on $\partial B_r(x_0)\setminus B_{\k r}(x_0)$, we have $w>\eps$ in $B_r(x_0)\setminus B_{\k r}(x_0)$. In particular $w\ge v_{\eps}$ on $\partial D_{\eps}$ where 
$$D_{\eps}:= (B_r(x_0)\setminus B_{\k r}(x_0))\cap \{u>\eps\}.$$
By comparison principle, we know $w\ge v_{\eps}$ in $D_{\eps}$ and since $w=v_{\eps} = \eps$ on $\partial B_{\k r}(x_0) \cap \{u>\eps\}$, hence from \eqref{green estimates}
\be\label{eq 2.6}
|\nabla v_{\eps}|\le |\nabla w|\le C(\mu,N,\k) \gamma\;\;\mbox{on $\partial B_{\k r}(x_0) \cap \{u>\eps\}$}.
\ee
Given that $\dive (A(x)\nabla v_{\eps})=0$ in $D_{\eps}$, we have by divergence theorem and \eqref{eq 2.6}
\[
\begin{split}
\int_{D_{\eps}} (A(x)\nabla v_{\eps})\cdot \nabla (v_\eps-u)\,dx&= \int_{\partial B_{\k r}(x_0)\cap \{u>\eps\}} (u-v_{\eps})(A(x)\nabla v_{\eps})\cdot \nu(y)\, d\HH^{N-1}(y)\\
&\le {C(\mu)} \int_{\partial B_{\k r}(x_0) \cap \{ u>\eps \}} |u-\eps||\nabla v_{\eps}|\,d\HH^{N-1}(y)\\
&\le {C(\mu,N,\k) \gamma}\int_{\partial B_{\k r}(x_0) \cap \{ u>\eps \}} |u-\eps|\,d\HH^{N-1}(y)=: M_0(u)
\end{split}
\]
justification of use of divergence theorem in $D_{\eps}$ can be found in \cite[equation (3.4)]{ACF84}. From the calculations above, we can write 

\begin{align*}
\;\; \;\;& \int_{D_{\eps}} (A(x)\nabla v_{\eps})\cdot \nabla (v_\eps-u)\,dx \le M_0\\
 \Rightarrow&  \int_{D_{\eps}} (A(x)\nabla v_{\eps})\cdot \nabla v_\eps\,dx  \le M_0+\int_{D_{\eps}} (A(x)\nabla v_{\eps})\cdot \nabla u\,dx\\
\Rightarrow &\mu \int_{D_{\eps}}|\nabla v_{\eps}|^2\,dx \le M_0 + \frac{1}{\mu} \int_{D_{\eps}} |\nabla v_{\eps}||\nabla u|\,dx\\
\Rightarrow & \mu \int_{D_{\eps}}|\nabla v_{\eps}|^2\,dx\le M_0+\frac{\eps_0}{2\mu} \int_{D_{\eps}} |\nabla v_{\eps}|^2\,dx+\frac{1}{2\eps_0 \mu}\int_{D_\eps} |\nabla u|^2\,dx
\end{align*}
putting very small $\eps_0>0$ in the last inequality, we have 
$$
\int_{D_{\eps}}|\nabla v_{\eps}|^2\,dx \le M_0+C(\mu)\int_{D_\eps} |\nabla u|^2\,dx=:M_1(u).
$$
Since $v_{\eps}=\eps \mbox{ in $B_{\k r}(x_0)\cap \{u>\eps\}$}$ which implies $\nabla v_{\eps}=0$ in $B_{\k r}(x_0)\cap \{u>\eps\}$ and $v_\eps = u $ in $B_{r}(x_0) \setminus D_{\eps}$, 
\be\label{b1}
\int_{B_r(x_0)}|\nabla v_\eps|^2\,dx = \int_{D_{\eps}}|\nabla v_{\eps}|^2\,dx+\int_{B_r(x_0)\setminus D_{\eps}}|\nabla u|^2\,dx=:M_2(u).
\ee
By the definition of $v_{\eps}$, $0< v_{\eps}\le u$ on $\partial D_{\eps}$ and $\dive(A(x)\nabla v_{\eps}) = \dive(A(x)\nabla u) = 0$ in $D_{\eps}$, therefore by comparison principle $0< v_{\eps} < u $ in $D_{\eps}$. In the set $B_{\k r}(x_0) \cap \{u>\eps\} $ we have $v_{\eps} =\eps<u$ and $v_{\eps} = u$ in $B_r(x_0)\cap \{u\le \eps\}$.  Overall we have $0< v_{\eps} \le u$ in $B_r(x_0)\cap \{u>\eps\}$ Therefore 
\be\label{b2}
\int_{B_r(x_0)}| v_\eps|^2\,dx  \le \int_{B_{r}(x_0) \cap \{u\le \eps\}} |u|^2\,dx + \int_{B_{r}(x_0) \cap \{u> \eps\}} |u|^2\,dx  = \int_{B_{r}(x_0)} | u|^2\,dx.
\ee
Hence, from \eqref{b1} and \eqref{b2}, $v_{\eps}$ is bounded in $H^1(B_r(x_0))$.  Therefore, up to a subsequence, there exists a limit $v = \lim_{\eps \to 0} v_{\eps} $ in weak $H^1$ sense, such that $v$ satisfies the following
\be\label{prop of v}
\begin{cases}
\dive(A(x)\nabla v) =0 \qquad \mbox{in $(B_r(x_0)\setminus B_{\k r}(x_0))\cap \{u>0\}$}\\
v = u\qquad \qquad \;\; \qquad \qquad \mbox{in $B_{r}(x_0)\cap \{u\le 0 \}$}\\
v =0 \qquad \;\; \qquad \qquad \qquad\mbox{in $B_{\k r}(x_0)\cap \{u>0\}$}\\
v=u \qquad \; \; \qquad \qquad \qquad \mbox{on $\partial B_r(x_0)$}.
\end{cases}
\ee
We verify the above properties \eqref{prop of v} of $v$ at the end of this proof.

Let us use the function $v$ as a test function with respect to minimality condition on $u$ in $B_r(x_0)$, we have
\[
\begin{split}
\int_{B_r(x_0)} \Big ( \< A(x)\nabla u,\nabla u\>+  \lambda (u)\Big ) \,dx& \le \int_{B_r(x_0)}\Big (  \< A(x)\nabla v,\nabla v\> + \lambda (v) \Big )\,dx\\
\end{split}
\]
since, $v=u$ in $\{u\le 0\}$ and $\{ v>0 \} \subset \{u>0\}$, the integration in the set $\{u\le 0\}$ gets cancelled from both sides and we are left with the terms mentioned below. \newline Set $D_0:= (B_r(x_0)\setminus B_{\k r}(x_0))\cap \{u>0\}$, we have
\[
\begin{split}
\int_{B_r(x_0)\cap \{u> 0\}} \Big (  \< A(x) \nabla u,\nabla u \> - \< A(x) \nabla v,\nabla v \> \Big )\,dx &\le \int_{B_{ r}(x_0)\cap \{u>0\}} (\Lambda(v)-\Lambda(u))\,dx\\
&= \int_{B_{\k r(x_0)}\cap \{u>0\}} (\Lambda(v)-\Lambda(u))\,dx\\ 
&=  \lambda_0 |B_{\k r}(x_0)\cap \{u>0\}| . \qquad(\lambda_0: = -(\lambda_ + - \lambda_-)). \\ 
\end{split}
\]
We have second equality above because $\chi_{\{v>0\}}= \chi_{\{u>0\}}$ in $D_0$. Since $v=0$ in $D_0$, we have 
\[
\begin{split}
\int_{B_{\k r}(x_0)\cap \{u>0\}} \< A(x) \nabla u,\nabla u \>\,dx + \int_{D_0} \Big ( \< A(x) \nabla u,\nabla u \> - \< A(x) \nabla v,\nabla v \> \Big )\,dx\le \lambda_0 |B_{\k r}(x_0)\cap \{u>0\}|.
\end{split}
\]
Using the ellipticity of $A$ and shuffling the terms in the above equation, we obtain
\be\label{eq2.5}
\begin{split}
\int_{B_{\k r}(x_0)\cap \{u>0\}} \Big (\mu |\nabla u|^2-\lambda_0 \Big ) \,dx &\le \int_{D_0} \Big ( \< A(x) \nabla v,\nabla v \>- \< A(x) \nabla u,\nabla u \> \Big ) \,dx\\
&= \int_{D_0} \Big ( \< A(x) \nabla (v-u), \nabla (v+u) \> \Big )\,dx\\
&=  \int_{D_0}  \Big (\< A(x) \nabla (v-u), \nabla (u-v+2v) \> \Big ) \,dx\\
& \le 2\int_{D_0} \< A(x)\nabla v, \nabla (v-u) \>\,dx \\
& \le  \liminf_{\eps\to 0} 2\int_{D_0} \< A(x) \nabla v_{\eps}, \nabla (v_{\eps}-u) \>\,\\
&= \liminf_{\eps\to 0} 2 \int_{D_{\eps}} \< A(x) \nabla v_{\eps}, \nabla (v_{\eps}-u) \>\,dx\qquad \mbox{(since $v_{\eps}=u$ in $D_{\eps}\setminus D_0$)}\\
&=\liminf_{\eps\to 0} 2\int_{\partial B_{\k r}(x_0)\cap \{u>\eps\}} (u-\eps) ( A(x)\nabla v_{\eps})\cdot \nu \,dx\\
&\le \liminf_{\eps\to 0} \frac{2}{\mu} \int_{\partial B_{\k r}(x_0)\cap \{u>\eps\}}(u-\eps)\big | \nu\cdot \nabla v_{\eps} \big |\,dx\,:=\,M.
\end{split}
\ee
The second to last equality in above calculation is obtained from integration by parts, its justification can be found in \cite[equation (3.4)]{ACF84}. From \eqref{eq2.5} and \eqref{eq 2.6},  and using the trace inequality in $H^1(B_{\k r})$ we have (for some different constant $C(\k)$),
\be\label{eq2.8}
\begin{split}
M &\le C(\mu,N,\k) \gamma \int_{\partial B_{\k r(x_0)}} u^+\,d\HH^{N-1}(x)\\
&\le C(\mu,N, \k) \gamma   \int_{B_{\k r(x_0)}} \Big (  |\nabla u^+|+\frac{1}{r} u^+ \Big )\,dx \\
& \le C(\mu,N, \k) \gamma \Bigg [ |B_{\k r(x_0)}\cap \{u>0\}|^{1/2}\left ( \int_{B_{\k r(x_0)}}|\nabla u^+|^2\,dx \right )^{1/2}+ \frac{1}{r}\sup_{B_{\k r}(x_0)} (u^+)\big | \{B_{\k r}(x_0)\cap \{u>0\}\}\big | \Bigg ]\\
& \le C(\mu,N, \k) \gamma \Bigg [  \frac{1}{2 \sqrt{-\lambda_0} } \int_{B_{\k r}(x_0)\cap \{u>0\}} |\nabla u^+|^2\,dx +  2 \sqrt{-\lambda_0}  |B_{\k r(x_0)}\cap \{u>0\}| \\ 
&\qquad \qquad  \qquad \qquad  \qquad \qquad \qquad \qquad  \qquad \qquad   \qquad \qquad \qquad + \frac{1}{r} \sup_{B_{\k r}(x_0)}(u^+) \int_{B_{\k r}(x_0)\cap \{u>0\}}1\,dx \Bigg ] \\
&= \frac{C(\mu,N, \k)\gamma }{2 \sqrt{-\lambda_0} }\left ( \int_{B_{\k r}(x_0)\cap \{u>0\}} |\nabla u^+|^2 - \lambda_0  \,dx \right )+\frac{C(\mu,N,\k)\gamma}{\lambda_0  r}\sup_{B_{\k r}(x_0)}(u^+) \int_{B_{\k r}(x_0)\cap \{u>0\}}\lambda_0   \,dx
\end{split}
\ee
we have used Hölder's inequality and then Young's inequality above. From Lemma \ref{subharmonic}, $u^+$ is $A-$subharmonic in $B_r(x_0)$. If $G'$ is the Green's function for $L'(v)=\dive(A(x)\nabla v)$ in $B_r(x_0)$, then by comparison principle and Green's representation
$$
u^+(x)\le \int_{\partial B_r(x_0)} u^+(y)  \left  ( A(y)\nabla_y G'(x,y)\right )\cdot \nu_y \,d\HH^{N-1}(y)\;\; \forall x\in B_{\k r}(x_0).
$$
Since for all $y\in \partial B_r(x_0)$ and $x\in B_{\k r}(x_0)$, we have $\frac{1}{|x-y|^{N-1}}\le \frac{C(\k)}{r^{N-1}}$, then using the Green's function estimates c.f. \cite[Theorem 3.3 (v)]{GW82} we get 
\be\label{eq2.9}
\begin{split}
\sup_{B_{\k r}(x_0)}u^+&\le C(\mu) \int_{\partial B_r(x_0)} \frac{u^+(y)}{|x-y|^{N-1}}\,d\HH^{N-1}(y)\\
&\le C(\mu,\k,N) \fint_{\partial B_r(x_0)} u^+\,d\HH^{N-1}(y)= C(\mu, \k,N)\gamma r.
\end{split}
\ee
Use \eqref{eq2.5} and \eqref{eq2.9} in \eqref{eq2.8} and we have
\[
\begin{split}
\mu \int_{B_{\k r}(x_0)\cap \{u>0\}} \Big ( |\nabla u|^2-\lambda_0  \Big )\,dx & \le 
\frac{C(\mu,\k,N)\gamma }{2 \sqrt{-\lambda_0} }\int_{B_{\k r}(x_0)\cap \{u>0\}}\Big (  |\nabla u^+|^2 -\lambda_0  \Big )\,dx\\ 
&\qquad+\frac{C(\mu,\k,N)\gamma}{\lambda_0  r}\sup_{B_{\k r}(x_0)}(u^+) \int_{B_{\k r}(x_0)\cap \{u>0\}}\lambda_0 \,dx\\
&\le \frac{C(\mu,\k,N) \gamma}{\mu \sqrt{-\lambda_0} } \left (  1+ \frac{C(\k)\gamma}{\sqrt{-\lambda_0}} \right )\mu\int_{B_{\k r}(x_0)\cap \{u>0\}} \Big ( |\nabla u|^2-\lambda_0 \Big ) \,dx.
\end{split}
\]
If $\gamma$ is small enough, then 
$$
\int_{B_{\k r}(x_0)\cap \{u>0\}} \Big ( |\nabla u|^2-\lambda_0 \Big ) \,dx=0
$$
in particular $|\{u>0\}\cap B_{\k r(x_0)}|=0$, that is $u^+=0$ almost everywhere in $B_{\k r}(x_0)$.

It remains to verify the properties of $v$ in \eqref{prop of v}. Before looking at the proof, we observe that for a given  $\vf \in C_c^{\infty}(D_0)$, then there exists $\eps_0>0$ such that $\vf \in C_c^{\infty}(D_{\eps})$ for all $\eps<\eps_0$. Indeed, since $\supp(\vf)$ is a compact set, and $\bigcup_{\eps>0} D_{\eps}$ is a cover of $\supp(\vf)$, then for a finite set $\{\eps_1, ...,\eps_n\}$ we have $\supp (\vf) \subset \bigcup_{i=1}^{n} D_{\eps_i} \subset D_{\eps_{max}}$ where $\eps_{max} = \max (\eps_1,...,\eps_n)$. Therefore, $\vf \in C_c^{\infty}(D_{\eps})$ for all $\eps<\eps_{max}$.

Let us first verify that $\dive(A(x)\nabla v)=0 $ in $D_0$. For this let $\varphi \in C_c^{\infty}(D_0)$, then from continuity of $u$, there exists a $\eps_0>0$ such that $\supp (\varphi) \subset D_\eps$ for all $\eps<\eps_0$, also we have 
\be\label{weak form}
\int_{D_0}\langle A \nabla v, \nabla \varphi \rangle \,dx = \int_{\supp (\varphi)} \langle A \nabla v, \nabla \varphi \rangle \,dx
\ee
since $\supp (\varphi) \subset D_\eps$, from the definition of $v_{\eps}$ we have 
$$
\int_{\supp (\varphi)} \langle A \nabla v_{\eps}, \nabla \varphi \rangle \,dx=0
$$
and we know that $v$ is a weak limit of $v_{\eps}$ in $H^1({B_r(x_0)})$, therefore from \eqref{weak form} we have 
$$
\int_{D_0}\langle A \nabla v, \nabla \varphi \rangle \,dx = \int_{\supp (\varphi)} \langle A \nabla v, \nabla \varphi \rangle \,dx= \lim_{\eps\to 0} \int_{\supp (\varphi)} \langle  \nabla v_{\eps}, \nabla \varphi \rangle \,dx=0.
$$
Hence we show that $\dive(A(x)\nabla v)=0 $ in $D_0$. To show that $v=0$ in $B_{\k r}(x_0)\cap \{u>0\}$, we now take the function $\varphi \in C_c^{\infty}(B_{\k r}(x_0)\cap \{u>0\})$. From the same reasoning as above we know that there exists an $\eps_0>0$ such that $\supp(\varphi)\subset B_{\k r}(x_0)\cap \{u>\eps\}$ for all $\eps<\eps_0$. From the definition of $v_{\eps}$, we have 
$$
\int_{\supp (\varphi)} v_{\eps} \varphi\,dx=\int_{\{u>\eps\} \cap B_{\k r}(x_0)} v_{\eps} \varphi\,dx = \eps \int_{\{u>\eps\} \cap B_{\k r}(x_0)}\varphi\,dx
$$
and in limit $\eps\to 0$, from the above equation we have 
$$
\int_{\supp(\varphi)} v\varphi\,dx= \lim_{\eps\to 0} \int_{\supp (\varphi)} v_{\eps} \varphi\,dx=0
$$
and therefore $v=0$ a.e. in $B_{\k r}(x_0)\cap \{u>0\}$. To prove that $v=u$ in $\{u\le 0\}$, we observe that $\{u\le 0\}\subset \{ u\le \eps\}$, hence from the definition of $v_{\eps}$, have
$$
v_{\eps} = u \mbox{ in $\{u\le 0\}$}.
$$
since the weak limits maintain the equality (c.f. \cite[Lemma 3.14]{MS211}) the claim follows in the limit $\eps \to 0$. Apart from that, since $v_{\eps}=u$ on $\partial B_r(x_0)$ therefore from conservation of traces in weak convergence, it follows that $v=u$ on $\partial B_r(x_0)$.
This completes the proof of Proposition \ref{non degeneracy}.
\end{proof}

\begin{remark}
In the proposition above the constant is local in nature, this means, the value of the constant depends on the choice of compact set $K\subset \subset B_2^+$ where $x_0\in K$.
\end{remark}

\begin{lemma}
Let $u_0$ and $u_k$ be as in Theorem \ref{convergence}. Then, for a subsequence of $u_k$, for any $R>0$ we have 
\be\label{pw convergence}
\chi_{\{u_k>0\} \cap B_R^+} \to \chi_{\{u_0>0\} \cap B_R^+}\qquad \mbox{ a.e. in $B_R^+$}.
\ee
This in turn implies 
\be\label{L1 convergence}
\chi_{\{u_k>0\} \cap B_R^+} \to \chi_{\{u_0>0\} \cap B_R^+}\qquad \mbox{ in $L^1(B_R^+)$}.
\ee
\end{lemma}

\begin{proof}
From Lemma \ref{convergence}, we can consider a subsequence of $u_k$ such that $u_k\to u_0$ in $L^{\infty}(B_R^+)$.
Let $x\in B_R^+$. If $x\in \{u_0>0\}\cap B_R^+$ (or $\chi_{\{u_0>0\}\cap B_R^+}(x)=1$), then  $u_k(x) > \frac{u(x)}{2}>0$ (or $\chi_{\{u_k>0\}\cap B_R^+}(x)=1$) for $k$ sufficiently large. Thus we conclude that 
$$
\mbox{ $\chi_{\{u_k>0\}\cap B_R^+} (x) \to \chi_{\{u_0>0\}\cap B_R^+} (x)$ as $k\to \infty$ for all $x\in \{u_0>0\}\cap B_R^+$}. 
$$
If $x\in \{u_0\le 0\}^o \cap B_R^+$ (or $\chi_{\{u_0>0\}\cap B_R^+}(x)=0$), then there exists $\delta>0$ such that $B_{\delta}(x) \subset \{u_0\le 0\} \cap B_R^+$. Thus we have $\frac{1}{\delta} \fint_{\partial B_\delta(x)} u_0^+ \,d\HH^{N-1} =0$. Again, by the uniform covergence of $u_k$ to $u_0$ in $B_R^+$ (c.f. Lemma \ref{convergence}) we obtain 
\be\label{c+}
\frac{1}{\delta} \fint_{\partial B_\delta(x)} u_k^+ \,d\HH^{N-1} \le \frac{1}{2}c(\mu,N,\lambda_{\pm})\qquad\mbox{for $k$ sufficiently large.}
\ee
Here $c(\mu,N,\lambda_{\pm})$ is as in Proposition \ref{non degeneracy}. This implies $u_k \le 0$ in $B_{\frac{\delta}{2}}(x)$ (c.f. Proposition \ref{non degeneracy}). In particular, $\chi_{\{u_k(x)\le 0\}}(x) =0$ for $k$ sufficiently large. This way, we obtain 
\be\label{c-}
\mbox{ $\chi_{\{u_k>0\}\cap B_R^+} (x) \to \chi_{\{u_0> 0\}\cap B_R^+} (x)$ as $k\to \infty$ for all $x\in \{u_0\le 0\}\cap B_R^+$}. 
\ee
From the representation theorem \cite[Theorem 7.3]{ACF84}, we know that $|\partial \{u_0>0\}\cap B_R^+| =0$. From \eqref{c+}, \eqref{c-} and the fact that $|\partial \{u_0>0\}\cap B_R^+| =0$, we obtain the claim \eqref{pw convergence}. Since $|\chi_{\{u_k>0\}\cap B_R^+}|\le 1$, the claim \eqref{L1 convergence} follows from Lebesgue's dominated convergence theorem.
\end{proof}

\section{The main result}
We rephrase the notion of the tangential touch of the free boundary to the fixed boundary, which is equivalent to the tangential touch condition mentioned in statement of Theorem \ref{main result}.

In the proof of our main result, we will show that given $u\in \PP_1$, for every $\eps>0$ there exists $\rho_{\eps}>0$ such that 
\be
\partial \{u>0\}\cap B_{\rho}^+ \subset B_{\rho}^+\setminus K_{\eps},\qquad \forall\; 0< \rho\le \rho_{\eps}
\ee 

where 
$$
K_{\eps}:= \big \{x \in \R^N_+ \,:\, x_N \ge \eps \sqrt{x_1^2+...+x_{N-1}^2}  \big \}.
$$

\begin{proof}[Proof of Theorem \ref{main result}]

We assume, by contradiction that the free boundaries of functions in $\PP_1$ do not touch the origin in a tangential fashion to the plane. Then, there exists $\eps>0$ and sequences $v_j \in \PP_1$ and $x_j\in F(v_j)\cap K_{\eps}$ such that 
$ |x_j|\to 0$ as $j\to \infty$. Let $r_j=|x_j|$ and we consider the blowups $u_j := (v_j)_{r_j}$.  

Let $u_0 := \lim_{j\to \infty}u_j$ as in Lemma \ref{convergence}.  Also, let $x_0\in \partial B_1^+ \cap K_{\eps}$ be a limit up to a subsequence (still called $x_j$) such that $x_0= \lim_{j\to \infty}\frac{x_j}{|x_j|}$. Since $x_j\in F(v_j)$, we have $v_j(x_j)=0$. Therefore on rescaling, $u_j(\frac{x_j}{r_j}) =\frac{1}{r_j}v_j(x_j)=0$.  In the limit as $j\to \infty$ we have
$$
u_0(x_0)=\lim_{j\to \infty} u_j \left (\frac{x_j}{r_j} \right )=0.
$$
From the density assumption that $u_j$ satisfy condition \eqref{density} and Lemma \ref{L1 convergence} we have for any given $R>0$
\be\label{density on rescaling}
\begin{split}
\frac{|\{u_0>0\}\cap B_R^+|}{|B_R^+|}= \fint_{B_R^+}\chi_{\{u_0>0\}}\,dx &=\lim_{j\to \infty}\fint_{B_R^+}\chi_{\{u_j>0\}}\,dx\\
&=  \lim_{j\to \infty} \frac{1}{|B_{R r_j}^+|} \int_{B_{Rr_j}^+}\chi_{\{v_j>0\}}\,dx\\
& =  \lim_{j\to \infty}\frac{|\{v_j>0\}\cap B_{R r_j}^+|}{|B_{R r_j}^+|} >\DD.
\end{split}
\ee
We can see that the computations done in \eqref{density on rescaling}, in fact shows that the density property remains invariant under blowup of any function $v$,. This way, we conclude that  the function $(u_0)_0$ which is the blowup limit of $(u_0)$ (in particular $(u_0)_0:= \lim_{r\to 0} (u_0)_r$) also satisfies 
\be\label{density 2}
\frac{ |\{ (u_0)_0>0 \}\cap B_R^+  |}{|B_R^+ |}>\DD,\qquad \forall R>0.
\ee

Now, we note that from Lemma \ref{convergence} $u_0\in \PP_{\infty}$. Moreover, from \eqref{density on rescaling} $u_0\not\equiv 0$ and from \cite[Theorem 4.2, Lemma 4.3]{KKS06} we have $u_0 \ge 0$ also, from \eqref{density 2}, we conclude $(u_0)_0 \not \equiv 0$. This way, again by \cite[Theorem 4.9]{KKS06}, we have $u_0(x)=c\, x_N^+$ for all $x\in \R^N_+$ for some constant $c>0$. 

Hence the function $u_0$ cannot be equal to zero at any point in $\R^N_+$. But we have $x_0\in \partial B_1^+\cap K_{\eps}$ and $u_0(x_0)=0$. This leads to a contradiction.
\end{proof}

\begin{remark}
Interested readers may check that the modulus of continuity $\sigma$ mentioned in the statement of Theorem \ref{main result} can be written as 
$$
\sigma(r) = \sup \left (  \frac{x_n}{\sqrt {x_1^2+...+x_{N-1}^2}} \,:\,x\in (B_{\rho}^+\setminus \{0\}) \cap F(u), \rho \le r \right ).
$$
\end{remark}

\newpage
 \newcommand{\noop}[1]{}

\end{document}